\documentclass[letterpaper, 10 pt, conference]{ieeeconf}  

\IEEEoverridecommandlockouts                              
\usepackage{epstopdf}
\usepackage{graphicx}
\usepackage{amsmath,amssymb,amsfonts}
\usepackage{algorithmic}
\usepackage{textcomp}
\usepackage{amstext}
\usepackage{xcolor,mathdots}
\usepackage{hhline}
\usepackage{multirow}
\usepackage{makecell}
\usepackage{float} 
\usepackage{subfigure}
\usepackage{url}
\usepackage{psfrag}

\newtheorem{proposition}{Proposition}

\title{\LARGE \bf
	A parameterized solution to the simultaneous stabilization problem
}

\author{Yufang Cui$^{1}$,  and Anders Lindquist$^{2}$
	\thanks{$^{1}$Department of Automation, Shanghai
		Jiao Tong University, Shanghai, China. {\tt\small cui-yufang@sjtu.edu.cn}}%
	\thanks{$^{2}$Department of Automation and School of Mathematical Sciences, Shanghai
		Jiao Tong University, Shanghai, China. {\tt\small alq@math.kth.se}}%
}

\begin{document}

\maketitle
\thispagestyle{empty}
\pagestyle{empty}

\begin{abstract}
In a series of fundamental papers BK Ghosh reduced the simultaneous stabilization problem to a Nevanlinna-Pick interpolation problem. In this paper we generalize some of these results allowing for derivative constraints. Moreover, we apply a method based on a Riccati-type matrix equation, called the Covariance Extension Equation, which provides a parameterization of all solutions in terms of a monic Schur polynomial. The procedure is illustrated by examples.
\end{abstract}

\section{Introduction}\label{Introduction}
Simultaneous stabilization is the problem of finding a single controller that stabilizes multiple plants \cite{Vidyasagar,Youla}. In this paper, we consider the following problem. Given a family $p_\lambda(s)$ of single-input single-output proper transfer functions of degree $n_\lambda$, represented as 
\begin{equation}\label{systems}
	p_\lambda(s)=\frac{\lambda x_1(s)+(1-\lambda)x_0(s)}{\lambda y_1(s)+(1-\lambda)y_0(s)}
\end{equation}
where $\lambda\in[0,1]$, $x_0(s),x_1(s),y_0(s),y_1(s)\in H$,  find a proper
compensator $k(s)$ such that the closed-loop
systems $p_\lambda(s)(1+k(s)p_\lambda(s))^{-1}$ are stable for all $\lambda\in[0,1]$. Here $H$ is the ring of proper rational functions with real coefficients with poles
in the open left half plane $\mathbb{C^-}$.

In \cite{Ghosh,Ghosh2}, BK Ghosh studied the simultaneous partial pole placement problem, which finds application in the design of a compensator for a family of linear dynamical
systems. To solve this problem, he proposed an interpolation method, which provides a new viewpoint to solve such problems.

In the present paper we apply a more general interpolation strategy based on our previous work on a Riccati-type approach to analytic interpolation \cite{CLccdc,CLtac}, which in turn is based on algorithms for the partial stochastic realization problem \cite{b2,BGuL,b1,BLpartial} and  on \cite{b6}. This allows for more general class of systems as there is a parameterization of all solutions with a prescribed interpolation points and derivative constraints are allowed.

More precisely, we  transform the simultaneous stabilization problem to an analytic interpolation problem. which in its most general (scalar) form can be formulated in the following way.   Given $m+1$ distinct complex numbers $z_0,z_1,\dots,z_m$ in the open unit disc $\mathbb{D}:=\{z \mid |z|<1\}$, consider the problem  to find a real Carath\'eodory function mapping the unit disc $\mathbb{D}$ to the open right half-plane, i.e., a real function $f$ that is analytic in $\mathbb{D}$ and satisfies $\text{Re}\{f(z)\}>0$ there, and which in addition satisfies the interpolation conditions 
\begin{align}
	\label{interpolation}
	\frac{f^{(k)}(z_{j})}{k!}=w_{jk},\quad&j=0,1,\cdots,m,   \\
	&   k=0,\cdots n_{j}-1 \notag
\end{align}
where $f^{(k)}$ is the $k$:th derivative of $f$, and the interpolation values $\{w_{jk}; j=0,1,\cdots,m, k=0,\cdots n_{j}-1\}$ are complex numbers that occur in conjugate pairs.  In addition we impose the complexity constraint that the interpolant  $f$ is rational of degree at most 
\begin{equation}
	\label{deg(f)}
	n:=\sum_{j=0}^{m}n_j -1.
\end{equation}
In general there are infinitely many solutions to this problem, but, as we shall see in Section~\ref{sec:CEE}, they can be completely parameterized in terms of an arbitrary $n$-dimensional Schur polynomial $\sigma(z)$. Freely choosing $\sigma(z)$ allows us to tune the solution to specifications.

The paper is organized as follows. In Section~\ref{sec:3condition}, we present necessary and sufficient conditions for a family of plants to be simultaneous stabilizable.  Section~\ref{sec:interpolation} shows how to transform the simultaneous stabilization problem to an analytic interpolation problem.  Section~\ref{sec:CEE} presents how to solve analytic interpolation problem based on the Covariance Extension Equation. In Section~\ref{sec:applications}, finally, we apply our method to some problems in simultaneous stabilization.

\section{Simultaneous stabilization problem}\label{sec:3condition}

To solve this problem we first collect some results based on the work of Ghosh \cite{Ghosh,Ghosh2}. To this end we first consider a special case: 
Given a pair of distict plants represented by coprime factorizations
\begin{equation}
\label{p0p1}
p_0(s)=\frac{x_0(s)}{y_0(s)}, \qquad p_1(s)=\frac{x_1(s)}{y_1(s)},
\end{equation}
where $x_{i}(s), y_{i}(s) \in H $ and $y_{i}(s)$ is proper but not strictly proper, find a proper compensator which can stabilize $p_0$ and $p_1$ simultaneously. Let $J$ be set of multiplicative units in $H$. That is, an element $u$ of $H$ is a multiplicative unit if there exists $v$ in $H$ such that $vu=uv=1$. Moreover, let $\mathbb{C^+}$ be closed right half of the complex plane including infinity.

\begin{proposition}
	 The pair of distinct plants $p_0,p_1$ is simultaneously stabilized by a proper compensator if and only if there exists $\Delta_{0}(s), \Delta_{1}(s) \in J$, such that the following holds.
	
\quad	(i) If $s_{1},s_2$, $\cdots, s_{t}$ are the zeros of $x_{0} y_{1}-x_{1} y_{0}$ in $\mathbb{C^+}$ with multiplicities $m_{1}, \cdots, m_{t}$, respectively, then $s_{1},s_2$, $\cdots, s_{t}$ must be the zeros of $\Delta_{0} y_{1}-\Delta_{1} y_{0}$ and $\Delta_{1} x_{0}-\Delta_{0} x_{1}$ with multiplicities at least $m_{1}, m_{2}, \cdots, m_{t}$, respectively.
	
\quad	(ii) If $x_{0} y_{1}-x_{1} y_{0}=0$ at $\infty$ with multiplicity $m_{\infty}$, then $\Delta_{1} x_{0}-\Delta_{0} x_{1}=0$ at $\infty$ with multiplicity $m_{\infty}$.
\end{proposition}

\medskip

\begin{proof}
	The main idea of proposition 1 is due to BK Ghosh \cite{Ghosh,Ghosh2}, and the proof that we now sketch is an adaptation of his procedure.
Let the required proper compensator be represented by the coprime factorization
\begin{equation}
\label{k}
k(s)=\frac{x_c(s)}{y_c(s)},
\end{equation}
 where $x_{c}(s),y_{c}(s)\in H$ and $y_{c}(s)$ is proper but not strictly proper. Then the transfer function of the closed-loop system is 
	\begin{equation}
		G_{i}(s)=\frac{n_{i}(s)}{d_{i}(s)}=\frac{x_i(s)y_c(s)}{y_i(s)y_c(s)+x_i(s)x_c(s)},\quad i={0,1}
	\end{equation}

Since $x_i(s),y_i(s),x_c(s),y_c(s)\in H$, which means all their poles are in $\mathbb{C^-}$, the poles of $n_{i}(s)$ and $d_{i}(s)$ are also in $\mathbb{C^-}$. To stabilize $p_0(s)$ and $p_1(s)$ simultaneously, we therefore need to have all zeros of $d_{i}(s)$ in $\mathbb{C^-}$.

If $x_i(s)/y_i(s)$ is proper but not strictly proper, then $n_{i}(s)$ is a nonzero number at infinity. To make $G_{i}(s)$ a proper function, $d_{i}(s)$ also needs to be a nonzero number at infinity. 

If $x_i(s)/y_i(s)$ is strictly proper, then $d_{i}(s)$ will be a nonzero number at infinity, and $G_{i}(s)$ will be strictly proper.

Therefore stabilizing $p_0$ and $p_1$ simultaneously relies on the existence of $\Delta_{0}, \Delta_{1} \in J$ such that
   \begin{equation}\label{cond1}
   	x_{i}(s) x_{c}(s)+y_{i}(s) y_{c}(s)=\Delta_{i}(s),\quad i={0,1}
   \end{equation}		
   Solving \eqref{cond1} for $x_{c}$ and $y_{c}$, we have
	\begin{equation}
		\label{8}
		\begin{aligned}
			&x_{c}(s)=(\Delta_{0} y_{1}-\Delta_{1} y_{0})/(x_{0} y_{1}-x_{1} y_{0}) \\
			&y_{c}(s)=(\Delta_{1} x_{0}-\Delta_{0} x_{1}) /(x_{0} y_{1}-x_{1} y_{0})
		\end{aligned}
	\end{equation}
    
	Condition (i) is necessary and sufficient for $x_{c}(s), y_{c}(s)$ to belong to $H$. Condition (ii) is necessary and sufficient for $y_c(s)$ to be proper but not strictly proper and $x_c(\infty)/y_c(\infty)\neq\infty$, which means $x_{c}(s) / y_{c}(s)$ is a proper rational function.
\end{proof}
\medskip
Next, we consider the more general case. Let \eqref{p0p1} be a pair of distinct plants. Consider
\begin{equation}
	p_\lambda(s)=\frac{x_\lambda(s)}{y_\lambda(s)}=\frac{\lambda x_1(s)+(1-\lambda)x_0(s)}{\lambda y_1(s)+(1-\lambda)y_0(s)}
\end{equation}
where $\lambda \in[0,1]$. Set $\eta_{ij}(s):=x_iy_j(s)-x_jy_i(s),i,j\in[0,1]$. 

\begin{proposition}
	The family of plants $p_\lambda(s)$ for $\lambda\in[0,1]$ is simultaneously stabilizable by a proper compensator if and only if there exists $\Delta_{0}$, $\Delta_{1}, \in J$ such that the conditions (i) and (ii) in Proposition 1 are satisfied together with the following additional condition:
\end{proposition}

\quad (iii) $\frac{\Delta_{1}}{\Delta_{0}}$ does not intersect the nonpositive real axis including infinity at any point in $\mathbb{C^+}$.

\medskip

\begin{proof}
	Let \eqref{k} be the required compensator. A necessary and sufficient condition for this compensator to stabilize the plants \eqref{p0p1} simultaneously is given by the conditions (i) and (ii). Additionally, \eqref{k} simultaneously stabilizes every other plant $x_{\lambda} / y_{\lambda}$ if and only if there exist $\Delta_{\lambda}$ $\in J, \lambda \in (0,1)$ such that
	\begin{equation}\label{cond2}
		x_{c} x_{\lambda}+y_{c} y_{\lambda}=\Delta_{\lambda}.
	\end{equation}
 By combining \eqref{cond1} and \eqref{cond2} we obtain
 \begin{equation}
 	 \lambda\Delta_1+(1-\lambda)\Delta_{0}=\Delta_{\lambda},\lambda\in(0,1)
 \end{equation}
 Since the poles of $\Delta_{0}$ and $\Delta_{1}$ are in $\mathbb{C^-}$, the poles of $\Delta_{\lambda}$ are in $\mathbb{C^-}$ as well. In order to have $\Delta_\lambda$ in $J$, we need that the zeros are in $\mathbb{C^-}$. From
 \begin{equation}
 	\lambda\Delta_1+(1-\lambda)\Delta_{0}=\Delta_{\lambda}=0
 \end{equation}
 we get 
 \begin{equation}
 	\frac{\Delta_{1}}{\Delta_{0}}=1-\frac{1}{\lambda} \in (-\infty,0), \text{for~} \lambda \in (0,1)
 \end{equation}
 
 If $\frac{\Delta_{1}}{\Delta_{0}}$ intersects $(-\infty,0)$ at $\mathbb{C^+}$ (suppose at $\hat{s}$), then there is a $\lambda\in(0,1)$, such that 
 \begin{equation}
 	\frac{\Delta_{1}}{\Delta_{0}}(\hat{s})=1-\frac{1}{\lambda}
 \end{equation}
 and 
 \begin{equation}
 	\Delta_{\lambda}(\hat{s})=0
 \end{equation}
 which means that $\Delta_{\lambda}$ is not in $J$ which contradicts the result that all $\Delta_{\lambda}\in J$. It follows that a necessary and sufficient condition for the existence of $\Delta_{\lambda}$ $\in J, \lambda \in (0,1)$  is given by the condition (iii) described above. 

\end{proof}
\medskip
\section{Determining the interpolation conditions}\label{sec:interpolation}

 In this section, we reformulate the three condition (i)-(iii) above as interpolation conditions. 
\begin{proposition}
	Let $s_{j}$ be a zero of $x_{0} y_{1}-x_{1} y_{0}$ in $\mathbb{C^+}$ of multiplicity $n+1$, but not a zero of $y_0$ and $y_1$, or $x_0$ and $x_1$. Then condition (i) and condition (ii) are equivalent the $i$-th derivative of $\Delta_{1}(s) / \Delta_{0}(s)$ satisfying the interpolation constraint
	\begin{equation}
		\label{prop3}
		(\frac{\Delta_{1}}{\Delta_{0}})^{(i)}(s_j)=(\frac{y_{1}}{y_{0}})^{(i)}(s_j)
	\end{equation}
where $\Delta_{1}, \Delta_{0} \in J$, $i=0,\cdots,n$.
\end{proposition}
\medskip
\begin{proof}
	We shall need the Leibniz formula 
	\begin{equation}
		[u(x) v(x)]^{(n)}=\sum_{k=0}^n C_n^k u^{(n-k)}(x) v^{(k)}(x),
	\end{equation}
	where $(n)$ is the $n$-th derivative. Since $s_{j}$ is a zero of $x_{0} y_{1}-x_{1} y_{0}$ of multiplicity $n+1$, 
	\begin{equation}
		\label{xy}
		(x_{0} y_{1}-x_{1} y_{0})^{(i)}(s_j)=0,\quad i=0,1,\cdots n.
	\end{equation}

	By condition (i), we need to have
	\begin{equation}
		\label{derivative}
		(\Delta_{0}y_1-\Delta_{1}y_0)^{(i)}(s_j)=0,\quad i=0,1,\cdots n
	\end{equation}

	which is equivalent to
	\begin{equation}
		\label{equal}
		(y_1-\frac{\Delta_{1}}{\Delta_{0}}y_0)^{(i)}(s_j)=0,\quad i=0,1,\cdots n.
	\end{equation}
	
	By Leibniz formula, \eqref{equal} implies
	\begin{equation}
		(y_1)^{(i)}(s_j)-\sum_{k=0}^{i} C_i^k (\frac{\Delta_{1}}{\Delta_{0}})^{(i-k)}(s_j) y_0^{(k)}(s_j)=0
	\end{equation}
	for $i=0,1,\cdots,n$. Thus, if $n=0$,
	\begin{equation}
		\frac{\Delta_{1}}{\Delta_{0}}(s_j)=\frac{y_{1}}{y_{0}}(s_j),
	\end{equation}
	Suppose that, for $i=0,\cdots,n-1$ and multiplicity $n$, \eqref{prop3} holds. Then for multiplicity $n+1$, we need an additional constraint
	\begin{equation}
		(y_1-\frac{\Delta_{1}}{\Delta_{0}}y_0)^{(n)}(s_j)=0
	\end{equation}
	which is 
	\begin{equation}
		(y_1)^{(n)}-\sum_{k=1}^{n} C_n^k (\frac{y_{1}}{y_{0}})^{(n-k)} y_0^{(k)}-(\frac{\Delta_{1}}{\Delta_{0}})^{(n)}y_0=0
	\end{equation}
	at $s_j$. By Leibniz formula,
	\begin{equation}
		(y_1)^{(n)}-(y_1)^{(n)}+(\frac{y_1}{y_0})^{(n)}y_0-(\frac{\Delta_{1}}{\Delta_{0}})^{(n)}y_0=0
	\end{equation}
	at $s_j$, yielding
	\begin{equation}
		(\frac{\Delta_{1}}{\Delta_{0}})^{(n)}(s_j)=(\frac{y_1}{y_0})^{(n)}(s_j).
	\end{equation}
	Then, by mathematical induction, Proposition 3 follows.
	Similarly, \eqref{xy} means
	\begin{equation}
		(\frac{x_{1}}{x_{0}})^{(i)}(s_j)=(\frac{y_1}{y_0})^{(i)}(s_j)\quad i=0,1,\cdots,n,
	\end{equation}
	so
	\begin{equation}
		(\frac{\Delta_{1}}{\Delta_{0}})^{(i)}(s_j)=(\frac{y_1}{y_0})^{(i)}(s_j)=(\frac{x_{1}}{x_{0}})^{(i)}(s_j),
	\end{equation}
	for $i=0,1,\cdots,n$ which concludes the proof
\end{proof}
\medskip
If $s_{j}$ is a zero of $ y_{1}$ and $ y_{0}$ with certain multiplicity, or a zero of $ x_{1}$ and $ x_{0}$ with certain multiplicity, then we need adjust the interpolation conditions. For example, let $s_{j}$ be a zero of $ x_{0}$ and $x_{1}$, but not a zero of $y_0$ and $y_1$. In this case, we need that $\Delta_{1} / \Delta_{0}$ interpolates the pair of numbers $(s_{j},(y_{1} /  y_{0})(s_{j}))$.

Condition (iii) means that $\frac{\Delta_{1}}{\Delta_{0}}$ maps $\mathbb{C^+}$ to the complex plane excluding the nonpositive real axis. In other words, the map is $\mathbb{C^+}\rightarrow re^{i\theta}, r\in(0,\infty),\theta\in(-\pi,\pi)$. 

Next we formulate the relevant analytic interpolation problem. Denote
\begin{equation}
	F(s):=\sqrt{\frac{\Delta_{1}}{\Delta_{0}}}
\end{equation}
which maps $\mathbb{C^+}$ to the open right half plane, i.e., to $\sqrt{r}e^{i\theta}, r\in(0,\infty),\theta\in(-\pi/2,\pi/2)$. 

 Using the Möbius transformation $z=(1-s)(1+s)^{-1}$, which maps $\mathbb{C^+}$ into the interior of the unit disc, we set
 \begin{equation}
 	f(z):=F((1-z)(1+z)^{-1})
 \end{equation}

Then the problem is reduced to finding a Carath\'eodory function $f(z)$ that  satisfies interpolation constraints. This is an analytic interpolation problem. Once we have solved for $f(z)$, we can do the following transformations to get the compensator $k(s)$:
\begin{equation}
	F(s)=f((1-s)(1+s)^{-1})
\end{equation}
\begin{equation}
	k(s)=\frac{F^2x_0-x_1}{y_1-F^2x_1}
\end{equation}

\section{The analytic interpolation problem}\label{sec:CEE}
In this section, we show how to solve the analytic interpolation problem \eqref{interpolation} using the Covariance Extension Equation \cite{CLccdc,CLtac,CLcdc}. 
To simplify calculations, we  normalize the problem by setting $z_0=0$ and $ f(0)=\tfrac{1}{2} $, which can be achieved through a simple M{\"o}bius transformation. Since $f$ is a real function,  $f^{(k)}(\bar{z}_j)/ k!=\bar{w}_{jk}$ is an interpolation condition whenever $f^{(k}(z_j)/ k!=w_{jk}$ is.

If $f$ is a Carath\'eodory function, then
\begin{equation}
	\label{phi+}
	\phi_+(z):= f(z^{-1})
\end{equation}
is a {\em positive real\/} function. The problem is then reduced  to finding a rational positive real function 
\begin{equation}
	\label{covexpansion}
	\phi_+(z)=\tfrac{1}{2}+c_{1}z^{-1}+c_{2}z^{-2}+c_{3}z^{-3}+\cdots ,
\end{equation}
of degree at most $n$ which satisfies the interpolation constraints \eqref{interpolation}.

Since $\phi_+(z)$ is analytic in $\mathbb{D}^{C}$ and $\phi_+(\infty)=\tfrac12$, there is an expansion
\begin{equation}\label{f}
	\phi_+(z)=\frac{1}{2}+c_{1}z^{-1}+c_{2}z^{-2}+c_{3}z^{-3}+\cdots ,
\end{equation}
and, since $\phi_+(z)$ is positive real, 
\begin{equation}
	\Phi(z):=\phi_+(z)+\phi_+(z^{-1})=\sum_{k=-\infty}^{\infty}c_{k}z^{-k}>0 \quad z\in \mathbb{T},
\end{equation}
where $\mathbb{T}$ is the unit circle $\{ z=e^{i\theta}\mid 0\leq\theta <2\pi\}$. Hence $\Phi$ is a power spectral density, and therefore there is a minimum-phase spectral factor $v(z)$ such that 
\begin{equation}
	\label{spectralfactor}
	v(z)v(z^{-1})=\Phi(z).
\end{equation}
Clearly $\phi_+$ has a representation
\begin{equation}\label{ab2f}
	\phi_+(z)=\frac{1}{2}\frac{b(z)}{a(z)} 
\end{equation}
where
\begin{subequations}
	\begin{equation}
		a(z)=z^{n}+a_{1}z^{n-1}+\cdots+a_{n}\\
	\end{equation}
	\begin{equation}
		b(z)=z^{n}+b_{1}z^{n-1}+\cdots+b_{n}\\
	\end{equation}
\end{subequations}
are Schur polynomials, i.e., monic polynomials with all roots in the open unit disc $\mathbb{D}$. Consequently 
\begin{equation}
	\label{GG*}
	v(z)v(z^{-1})=\frac{1}{2}\left[\frac{b(z)}{a(z)} +\frac{b(z^{-1})}{a(z^{-1})}\right],
\end{equation}
and therefore 
\begin{equation}
	\label{G}
	v(z)=\rho\frac{\sigma(z)}{a(z)},
\end{equation}
where $\rho>0$ and 
\begin{equation}\label{sigma}
	\sigma(z)=z^{n}+\sigma_{1}z^{n-1}+\cdots+\sigma_{n} 
\end{equation}
is a Schur polynomial. It follows from \eqref{GG*} and \eqref{G} that
\begin{equation}
	\label{ab2sigma}
	a(z)b(z^{-1})+b(z)a(z^{-1})=2\rho^{2}\sigma(z)\sigma(z^{-1}).
\end{equation}
We shall represent the monic polynomials $a(z)$, $b(z)$ and $\sigma(z)$ by the $n$-vectors $a=[a_1,a_2,\cdots,a_n]',b=[b_1,b_2,\cdots,b_n]',\sigma=[\sigma_1,\sigma_2,\cdots,\sigma_{n}]'$.

Following \cite{b6} we note that \eqref{ab2f} has an observable realization
\begin{equation}
	\phi_+(z)=\frac{1}{2}+h'(zI-F)^{-1}g
\end{equation}
where  
\begin{subequations}
	\begin{equation}\label{g}
		F=J-ah',\quad  g=\frac{1}{2}(b-a),
	\end{equation}
	\begin{equation}\label{J}
		h=\begin{bmatrix}1\\0\\\vdots\\0\end{bmatrix},
		\quad
		J=\begin{bmatrix}
			0&1&0&\cdots&0\\
			0&0&1&\cdots&0\\
			\vdots&\vdots&\vdots&\ddots&\vdots\\
			0&0&0&\cdots&1\\
			0&0&0&\cdots&0\\
		\end{bmatrix}.
	\end{equation}
\end{subequations}

From stochastic realization theory  \cite[Chapter 6]{LPbook} it follows that the minimum-phase spectral factor \eqref{G} has a realization
\begin{equation}
	v(z)=\rho+h'(zI-F)^{-1}k
\end{equation}
where
\begin{equation}
	\rho=\sqrt{1-h'Ph},\quad k=\rho^{-1}(g-FPh)
\end{equation}
with $P$ being the minimum solution of the algebraic Riccati equation
\begin{equation}\label{ARE}
	P=FPF'+(g-FPh)(1-h'Ph)^{-1}(g-FPh)' .
\end{equation}
Following the calculations in \cite{BLpartial,b6} we now see that
\begin{equation}\label{gk}
	g=\Gamma Ph+\sigma-a,\quad k=\rho(\sigma-a) 
\end{equation}
and that \eqref{ARE} can be reformulated as 
\begin{equation}\label{P}
	P=\Gamma(P-Phh'P)\Gamma'+gg' 
\end{equation}
where $\Gamma$ is given by 
\begin{equation}
	\label{Gamma}
	\Gamma =J-\sigma h' .
\end{equation}

The  analytic interpolation problem amounts to finding $(a,b)$ given some interpolation data $w_{jk}$ and a particular Schur polynomial $\sigma(z)$.  

In \cite{CLtac}, we derive the condition for the existence of solutions of the analytic interpolation problem, which only depend on the interpolation data (see Proposition 5 in \cite{CLtac}). 

If the solution exists, \cite{CLtac} also shows that the {\em Covariance Extension Equation (CEE)}
\begin{subequations}\label{PgCCE}
	\begin{equation} \label{CEE}
		P = \Gamma (P-Phh'P) \Gamma' + g(P)g(P)' 
	\end{equation}
	(where $^\prime$ denotes transposition) with
	\begin{equation}\label{g(P)} 
		g(P)= u +U\sigma + U\Gamma Ph ,
	\end{equation}
\end{subequations}
where $u$ and $U$ are totally determined by the interpolaton data \eqref{interpolation}, has a unique symmeric solution $P\geq 0$ such that $h'Ph<1$. 
Moreover, for each $\sigma$ there is a unique solution of the analytic interpolation problem, and it is given by
\begin{subequations}\label{Psigma2ab}
	\begin{equation}\label{a}
		a=(I-U)(\Gamma Ph+\sigma)-u
	\end{equation}
	\begin{equation}\label{b}
		b =(I+U)(\Gamma Ph+\sigma)+u
	\end{equation}
	\begin{equation}\label{rho}
		\rho=\sqrt{1-h'Ph} ,
	\end{equation}
\end{subequations}
and 
the degree of $f(z)$ equals the rank of $P$. To solve equation \eqref{PgCCE}, a homotopy continuation method can be used, more details can be found in \cite{CLtac}. From above, we can easily draw the conclusion that if the interpolation data of the simultaneous stabilization problem satisfies the condition for the existence of the solution. Then different choices of Schur polynomial $\sigma(z)$ can generate different feasible solutions.

\section{ Computational examples}\label{sec:applications}
\subsection{Example 1}

\begin{figure}[htp]
	\centering
	\includegraphics[width=1\linewidth]{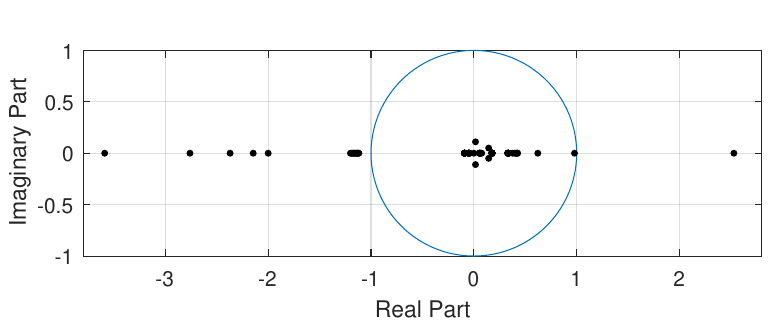}
	\caption{The poles of $p_\lambda$ before stabilization}
	\label{beforestabilization}
\end{figure}

Let us consider a simple case. Given $x_0,y_0,x_1,y_1\in H$ as 
\begin{equation}
	x_0=\frac{(s-15)(s-6)}{(s+0.5)(s+1.2)},\quad y_0=\frac{(s-3)(s-18)}{(s+1.5)(s+0.3)}
\end{equation}
\begin{equation}
	x_1=\frac{(s+9)(s-2)}{(s+0.7)(s+1.1)},\quad y_1=\frac{(s-11)(s+1)}{(s+0.9)(s+0.4)}
\end{equation}
 there are unstable poles when $\lambda$ varies on the interval $(0,1)$. To show the poles more clearly, we do the following transformation:
\begin{equation}\label{trans}
	z=\frac{1+s}{1-s}, 
\end{equation}
which maps the left half plane to the inside of the unit circle and maps the right half plane to the outside of the unit circle. Then a stable system has all poles inside the unit circle. After transformation \eqref{trans}, we can show  all poles of $p_\lambda$  when $\lambda$ varies from 0 to 1 at
intervals of 0.1 in Fig.~\ref{beforestabilization}. \\

From Fig.~\ref{beforestabilization}, we can see there are some systems that are not stable. Using the method in this paper, we first observe that $x_0y_1-x_1y_0$ has two zeros at $s_0=3169/165$ and $s_1=1113/250$ in $\mathbb{C^+}$. To make the systems stable, we therefore need the interpolation conditions
\begin{equation}
	(\frac{\Delta_{1}}{\Delta_{0}})(s_0)=(\frac{y_1}{y_0})(s_0),\quad(\frac{\Delta_{1}}{\Delta_{0}})(s_1)=(\frac{y_1}{y_0})(s_1)
\end{equation}
Using the Möbius transformation $z=(1-s)(1+s)^{-1}$, which maps  the open right half plane into the interior of the unit disc, the problem is reduced to finding a Carath\'eodory function  $f(z)$ that satisfies 
\begin{equation}
	f(\frac{1-s_0}{1+s_0})=\sqrt{(\frac{y_1}{y_0})(s_0)},\quad f(\frac{1-s_1}{1+s_1})=\sqrt{(\frac{y_1}{y_0})(s_1)}
\end{equation}
This is a Nevanlinna-Pick interpolation problem, which is a speical case of the analytic interpolation problem with $n_0=n_1=1$. By Proposition 5 in \cite{CLtac}, there exists solutions. Here we choose $\sigma(z)=z-0.9$. After calculation, we can get 

\begin{equation}
\frac{\Delta_1}{\Delta_{0}}=\frac{19.871 (s+0.1023)^2}{(s+9.988)^2}
\end{equation}

\begin{figure}[htp]
	\centering
	\includegraphics[width=0.8\linewidth]{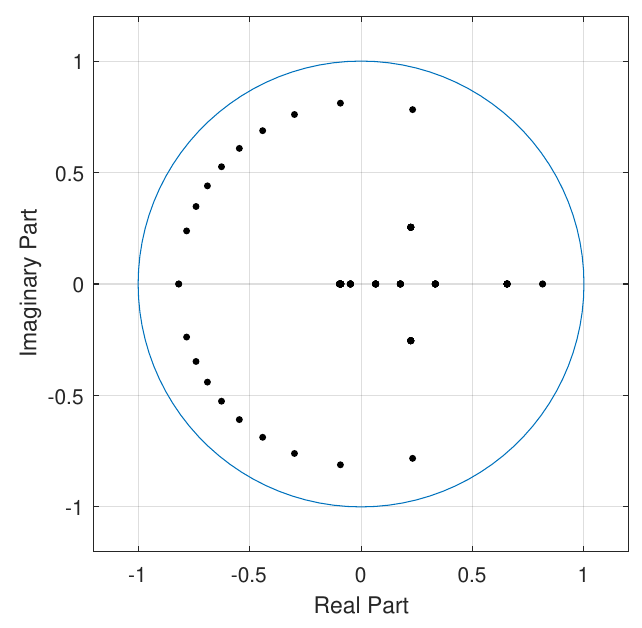}
	\caption{The poles of $p_\lambda$ after stabilization}
	\label{afterstabilization}
\end{figure}

After stabilization, the poles of $p_\lambda,\lambda\in[0,1]$ are shown in Fig.~\ref{afterstabilization}. Since all poles are in the open unit disc,  all feedback systems are stable.

To verify that different choices of $\sigma(z)$ produce different feasible solutions, let us vary the zero $z_0$ of $\sigma(z)$ from 0 to 1. Fig.~3 shows the results with $z_0=0,0.2,0.4,0.6,0.8,0.99$ respectively. We can see the solution changes with different $\sigma(z)$.
\begin{figure}[htbp]
	\centering
	\subfigure[$z_0=0$]{
		\begin{minipage}[t]{0.5\linewidth}
			\centering
			\includegraphics[width=1\linewidth]{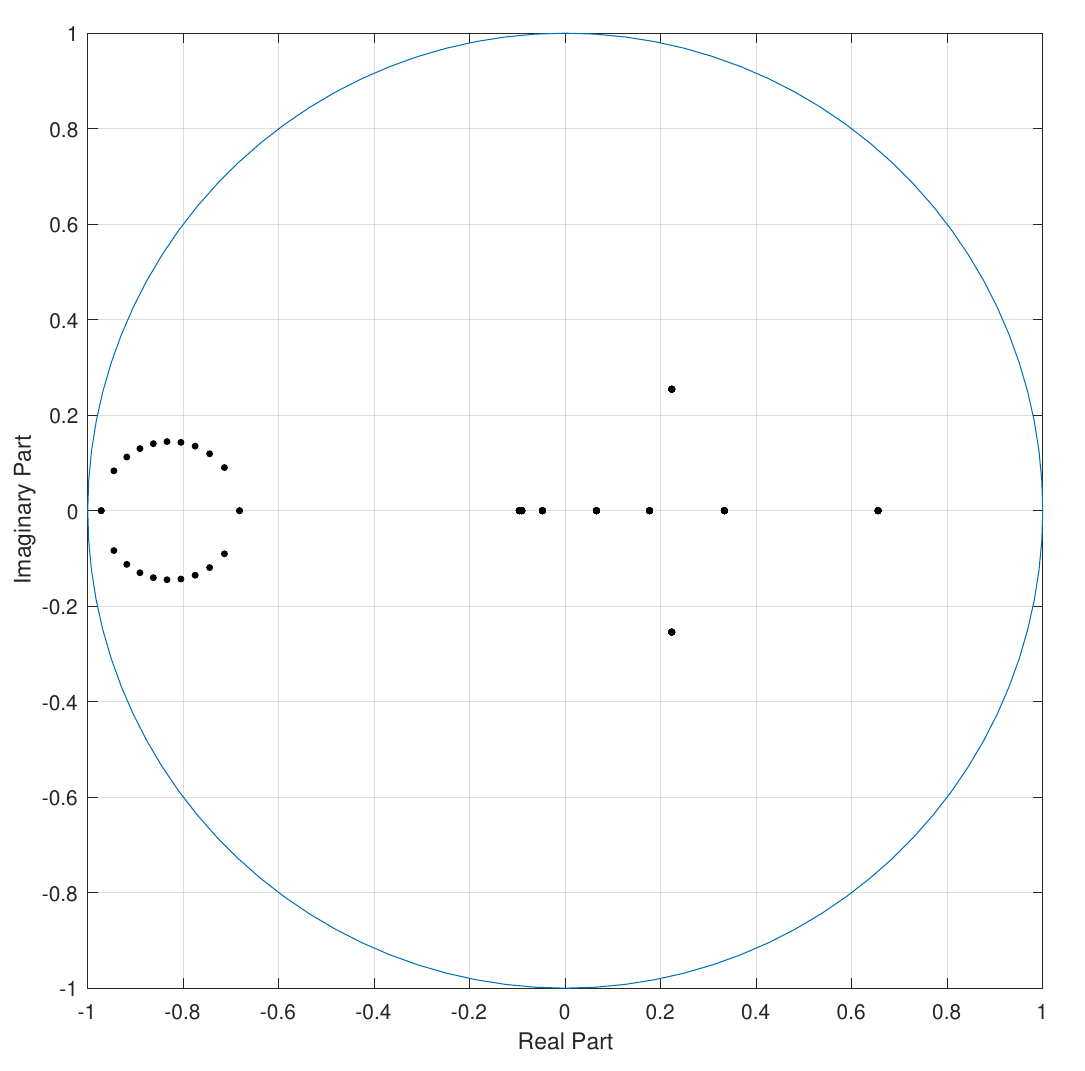}
		\end{minipage}%
	}%
	\subfigure[$z_0=0.2$]{
		\begin{minipage}[t]{0.5\linewidth}
			\centering
			\includegraphics[width=1\linewidth]{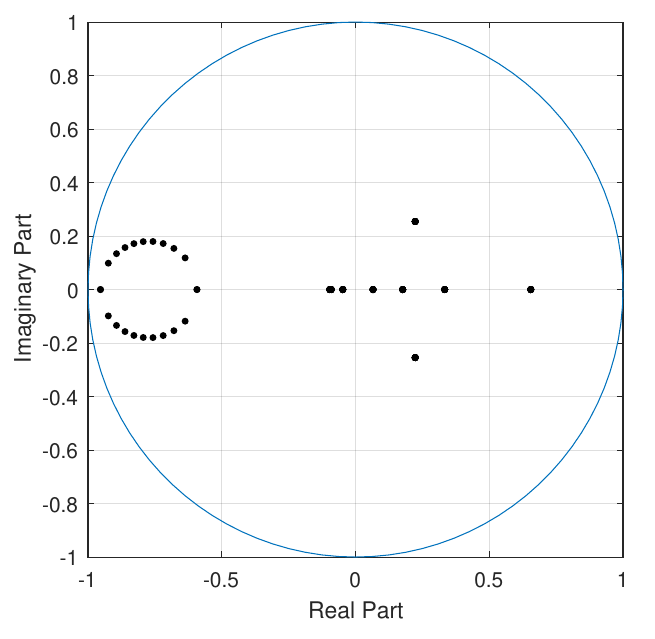}
		\end{minipage}%
	}%

	\subfigure[$z_0=0.4$]{
		\begin{minipage}[t]{0.5\linewidth}
			\centering
			\includegraphics[width=1\linewidth]{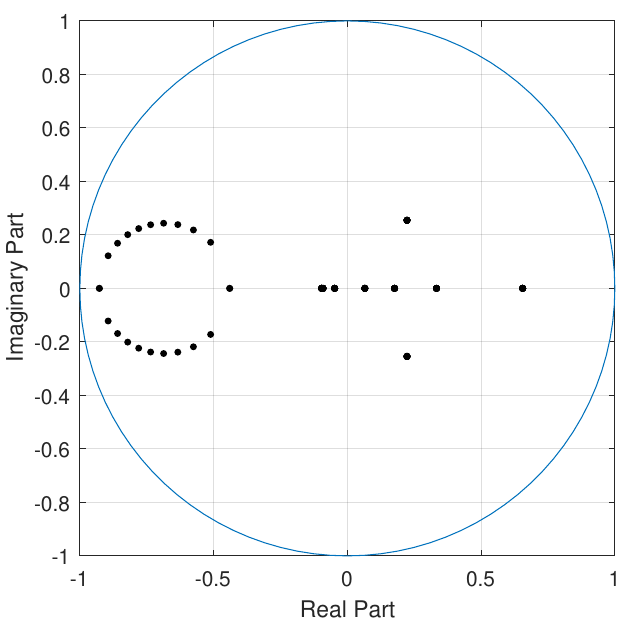}
		\end{minipage}
	}%
	\subfigure[$z_0=0.6$]{
		\begin{minipage}[t]{0.5\linewidth}
			\centering
			\includegraphics[width=1\linewidth]{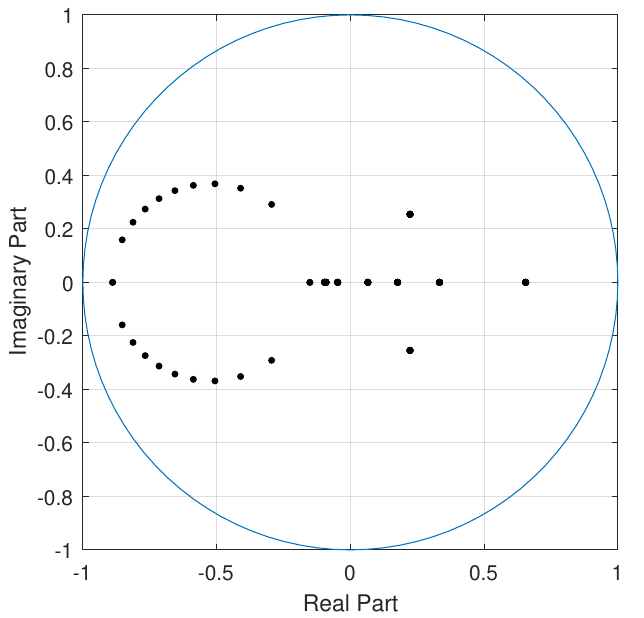}
		\end{minipage}
	}

		\subfigure[$z_0=0.8$]{
		\begin{minipage}[t]{0.5\linewidth}
			\centering
			\includegraphics[width=1\linewidth]{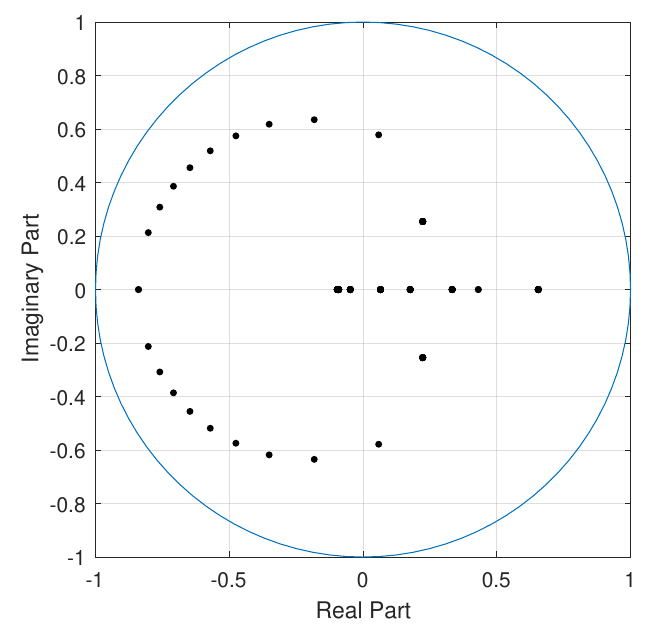}
		\end{minipage}
	}\subfigure[$z_0=0.99$]{
		\begin{minipage}[t]{0.5\linewidth}
			\centering
			\includegraphics[width=1\linewidth]{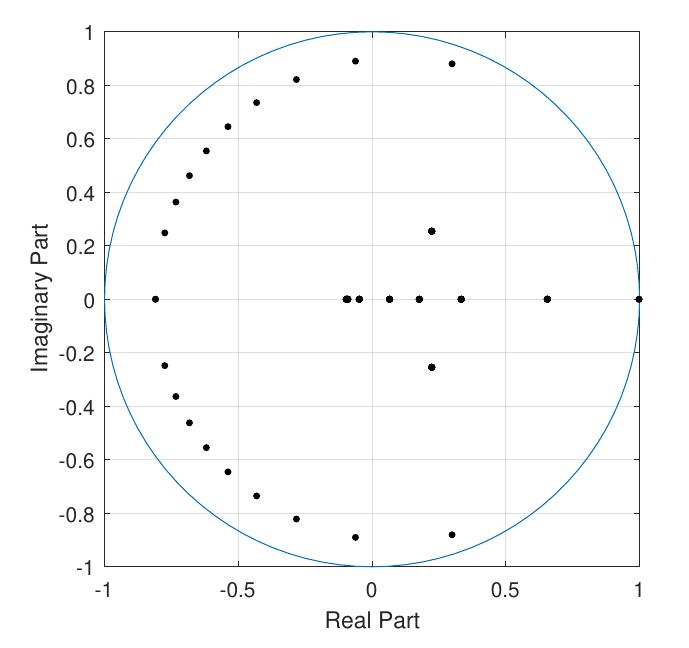}
		\end{minipage}
	}%
	\centering
	\label{sigma_zero}
	\caption{The poles of the stabilized system with diferent $\sigma(z)$}
\end{figure}

\subsection{Example 2}
Next we consider  more complex systems which includes derivative constraints, namely
\begin{equation}
	x_0=\frac{(s-0.2)(s+0.5)}{(s+0.3)(s+0.7)},\quad y_0=\frac{(s-1)^2}{(s+1.7)(s+0.2)}
\end{equation}
\begin{equation}
	x_1=\frac{2(s-0.2)(s+1.2)}{(s+0.4)(s+1.4)},\quad y_1=\frac{(s-1)^2}{(s+1.1)(s+0.6)}
\end{equation}
Obviously, when $\lambda$ varies from 0 to 1, $p_\lambda$ has poles at 1, which means all systems are unstable.

\begin{figure}[htp]
	\centering
	\includegraphics[width=0.8\linewidth]{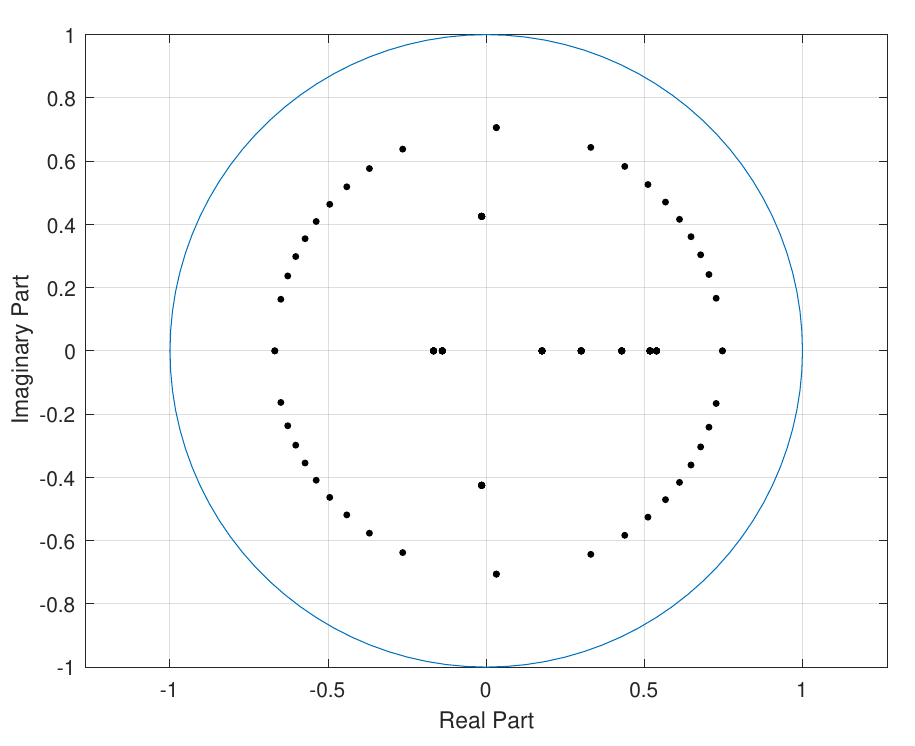}
	\caption{The poles after stabilization}
	\label{example2}
\end{figure}

By calculation, $x_0y_1-x_1y_0$ has zeros at 1 and 0.2 with multiplicity 2 and 1  respectively. This means that we need to find $\Delta_0$ and $\Delta_1$ such that
\begin{align}
	&\frac{\Delta_{1}}{\Delta_{0}}(1)=\frac{x_1}{x_0}(1),\quad(\frac{\Delta_{1}}{\Delta_{0}})'(1)=(\frac{x_1}{x_0})'(1)\\
	&\frac{\Delta_{1}}{\Delta_{0}}(0.2)=\frac{y_1}{y_0}(0.2)
\end{align}
Using the Möbius transformation $z=(1-s)(1+s)^{-1}$, which maps the open right half plane into the interior of the unit disc, the problem is reduced to finding a Carath\'eodory function  $f(z)$ that satisfies 
\begin{align}
	&f(0)=\sqrt{\frac{x_1}{x_0}(1)},\quad f'(0)=-(\frac{x_1}{x_0})'(1)/f(0)\\
	&f(\frac{2}{3})=\sqrt{\frac{y_1}{y_0}(0.2)}
\end{align}
which is an analytic interpolation problem with derivative constraint. Here we choose $\sigma(z)=z(z-0.1)$. By calculation, we can get 
\begin{equation}
	\frac{\Delta_1}{\Delta_0}=\frac{0.26463 (s+5.034)^2 (s+0.1448)^2}{(s^2 + 0.6404s + 0.9181)^2}
\end{equation}
Since there are three interpolation constraints,  we can get an $f(z)$ of degree 2 and a $\frac{\Delta_1}{\Delta_0}$ of degree 4. The poles of all $p_\lambda$ with $\lambda$ changing from 0 to 1 at interval 0.1 are showed in Fig.~\ref{example2}. Since all poles are in the open unit disc,  all systems are stable.

\section{conclusion}
In this paper we have studied the simultaneous stabilization problem to find a feedback compensator which stabilizes all the SISO systems \eqref{systems}. This leads to an analytic interpolation problem, which we solve by using a Riccati-type  algebraic matrix equation, the Covariance Extension Equation. This problem has infinitely many solutions, but we provide all of them, parameterized by a monic Schur polynomials. In future work we shall generalized these results to the MIMO case.
\bibliographystyle{IEEEtran}

\end{document}